\documentclass{kluwer}    
\usepackage{epsfig}
\usepackage{graphicx}
\usepackage{color}

\newdisplay{guess}{Conjecture}
\newdisplay{lemma}{Lemma}
\newtheorem{theorem}{Theorem}
\newtheorem{definition}{Definition}
\newproof{remark}{Remark}
\newproof{proof}{Proof}

\begin{document}
\begin{article}
\begin{opening}
\title{Symmetries in Differential Geometry: A Computational Approach to Prolongations}
\author{Thomas \surname{Branson}}
\author{Alfredo \surname{Villanueva}}
\runningauthor{Thomas Branson and Alfredo Villanueva}
\runningtitle{Symmetries in Differential Geometry: A Computational Approach to Prolongations}
\institute{The University of Iowa}
\date{May 3, 2007}

\begin{abstract}
The aim of this work is to develop a systematic manner to close overdetermined systems arising from conformal Killing tensors (CKT). The research performs this action for 1-tensor and 2-tensors. This research makes it possible to develop a new general method for any rank of CKT. This method can also be applied to other types of Killing equations, as well as to overdetermined systems constrained by some other conditions.

The major methodological apparatus of the research is a decomposition of the section bundles where the covariant derivatives of the CKT land via generalized gradients. This decomposition generates a tree in which each row represents a higher derivative. After using the conformal Killing equation, just a few components (branches) survive, which means that most of them can be expressed in terms of lower order terms. This results in a finite number of independent jets. Thus, any higher covariant derivative can be written in terms of these jets.

The findings of this work are significant methodologically and, more specifically, in the potential for the discovery of symmetries. First, this work has uncovered a new method that could be used to close overdetermined systems arising from conformal Killing tensors (CKT). Second, through an application of this method, this research finds higher symmetry operators of first and second degree, which are known by other means, for the Laplace operator. The findings also reveal the first order symmetry operators for the Yamabe case. Moreover, the research leads to conjectures about the second order symmetries of the Yamabe operator.
\end{abstract}


\end{opening}

\section{Introduction}

The search for a unified theory in physics is resulting in new concepts and theories. Some of these new developments have included overdetermined systems for partial differential equations, conformal Killing tensors and symmetries.

Overdetermined systems for ordinary differential equations (ODE) are very well known and they follow a standard procedure called "prolongations". This method is as follows: if $y^{n}=F(x,y,y^{'},...,y^{n-1})$, is an ODE, one introduces new variables, such that $y=z_{1} , z_{1}^{'}=z_{2},..., z_{n-1}^{'}= F(x,z_{1},...,z_{n-1}) $, reducing to first order differential equation. This also known as a closed system. However, in the case of overdetermined systems in partial differential equations, the method of prolongation is much more complex. As a result, usually ad hoc methods are used.

 Our research has developed a systematic manner for doing prolongation over certain types of overdetermined systems in partial differential equations. This research has applications finding hidden symmetries that agree with those found by other means in the Laplacian and Yamabe case. Notice that a symmetry of a differential operator L is a linear differential operator $D$ so that $LD=\hat{D}L$ for some linear differential operator $\hat{D}$ \cite{east1}. Furthermore, this method allows one to conjecture about symmetry operators of  second degree for the conformal Laplacian.

In \cite{b1} the authors show explicitly the prolongation for some examples, and predict the form of a prolongation for a large class of examples, without computing the prolongation in details. This method gives us a way to know a priori the number of new variables we require for closed systems.

To explicitly carry out prolongations in a systematic process, this research uses representation theory and generalized gradients, which allow us to express all higher derivatives in terms of the new variables (jets). These new variables land in sub-bundles arising from an irreducible representation of the group structure background.

Roughly speaking, we first express the general form of possible higher symmetry operators to the correspondent Laplace or Yamabe operator. Next, we express all these higher derivatives ($LD$, $\hat{D}L$) in terms of the new jets. Finally a comparison of both results ($LD-\hat{D}L=0$) yields the exact symmetry operators.

This research has used the Ricci program, a Mathematica package for doing tensor calculations in differential geometry, by John M. Lee \cite{lee}.

\section{Geometry}\label{s2}

We are working on a Riemannian n-manifold with metric $g_{ab}$, and $\nabla_{a}$ is the corresponding connection associated to g. We are using the Einstein summation convention; where repeated indices means a sum over these indices.

\subsection{Origin of the Yamabe Operator}
Re-scaling the metric in a Riemannian (or pseudo-Riemannian) manifold M of dimension n and metric g means
$$\widehat{g}=e^{2\rho}g,$$ where $ \rho$ $\in$ $C^{\infty}(M)$. In
other words, the original metric is multiply by a positive function. Such change of metrics are study in conformal geometry and mathematical physics, eg. \cite{east2} and \cite{fegan} . Important elements in differential geometry are;
$$
\begin{array}{l}  \\ \\ \texttt{Riemannian tensor} \\ \texttt{Ricci tensor} \\ \texttt{Scalar curvature} \end{array}
\begin{array}{c} \texttt{Metric $g$} \\ \\ {R}^{a}{}_{bcd} \\ {R}_{ab} \\ Sc \end{array}
\begin{array}{c}  \\ \\ \\ \\ \end{array}
\begin{array}{c} \texttt{Metric $\widehat g$} \\ \\ \widehat{R}^{a}{}_{bcd} \\\widehat{R}_{ab} \\ \widehat{Sc} \end{array}
$$
These elements might be related as,

$$ \widehat{R}^{a}{}_{bcd}=R^{a}{}_{bcd}+ \rho_{bc}\delta^{a}{}_{d}+ \rho_{bd}\delta^{a}{}_{c}-g_{bc}\rho^{a}{}_{d}+g_{bd}\rho^{a}{}_{c} $$
Where $\rho_{bc}=\rho_{b.c}-\rho_{b}\rho_{c}+\frac{1}{2}g^{ad}\rho_{a}\rho_{d}g_{bc} $ and $ \rho_{a}=\frac{\partial \rho}{\partial x^{a}}$.
$$ \widehat{R}_{bc}=R_{bc}-(n-2)\rho_{bc}-g_{bc}\rho^{a}{}_{a}.$$
$$ \widehat{Sc}=e^{-2\rho}(Sc-2(n-1)\rho^{a}{}_{a}).$$

The Laplace-Beltrami operator ($ \Delta$) in terms of the previous formulas can be written as \cite{yamabe}:
$$
\begin{array}{rl}
\widehat{Sc}&=e^{-2\rho}(Sc-2(n-1)\rho^{a}{}_{a})\\\\
&=e^{-2\rho}(Sc-2(n-1)[\Delta\rho+(\frac{n}{2}-1)g^{ab}\rho_{a}\rho_{b}])\\\\
&=e^{-2\rho}(Sc-\frac{4(n-1)}{n-2}e^{-(\frac{n}{2}+1)\rho}\Delta e^{(\frac{n}{2}-1)\rho}),
\end{array}
$$
and setting $ u=e^{(\frac{n}{2}-1)\rho}$
$$ \widehat{Sc}u^{\frac{n+2}{n-2}}= [Sc-\frac{4(n-1)}{n-2}\Delta]u$$
or
$$ \widehat{Sc}\frac{n-2}{4(n-1)}u^{\frac{n+2}{n-2}}= [\underbrace{-\Delta+\frac{n-2}{4(n-1)}Sc}]u.$$

The right side of this formula is know as the Yamabe operator (or Conformal Laplacian), and it is a conformal differential operator on Riemannian manifolds ($M^{n}$,g) for n$\geq $3, written in common literature as:
$$Yf=\Delta f+\frac{n-2}{4(n-1)}Scf,$$
where $\Delta $ is the Bochner Laplacian (negative of the Beltrami-Laplace operator). This last convention for $\Delta$ is used in the present paper.

\section{Associated Bundles and some Representation Theory} \label{ss1}

Associated bundles basically are principal G-bundles (G is a topological group) where each fiber G is changed by a fiber $V$. The group G is known as the group structure. In a Riemannian manifold (M,g), Oriented Riemannian manifold or in a Spin manifold the group structure G is O(n), SO(n) or Spin(n) respectively. In this work we are using O(n) group structure.

\subsection{Associated Bundles}
Let G be a topological group, $ \pi : B \longrightarrow M $ a principal G-bundle; where M is a Riemannian manifold, V a topological space and
$$\lambda : G \longrightarrow Aut(V),$$
a representation G. Then the fiber bundle associated to B by $\lambda$, is a fiber bundle $\pi_{\lambda}: B \times_{\lambda}V \longrightarrow M$ with fiber V and group G, that is defined as follows:
$$ B[\lambda] = B \times_{\lambda}V = B \times_{\lambda}\frac{V}{G}$$
where the action of G on $B \times V $ is defined by
$$ g.(p,v)=(pg^{-1},\lambda(g)v), \quad \forall g \in G, \quad p \in M, \quad v \in V $$
$\pi_{\lambda} $ is defined by
$$ \pi_{\lambda}[p,v] = \pi(p),$$
where [p,v] represent the G-orbit of (p,v) $\in B \times V.$

\subsection{Generalized Gradients}
Generalized gradients are first-order differential operators and their origin goes back to Stein and Weiss in \cite{sw}, Fegan in \cite{fegan} and T. Branson in \cite{b2}. There is a generalization to an arbitrary semi simple Lie group G by Orsted in \cite{orsted}. We use these results for G = O(n) in the following fashion:
Let M be a Riemannian manifold ($M^{n}$,g), consider the orthonormal frame bundle E. Let $B[\lambda] = B \times_{\lambda}V$ be the corresponding associated bundle. Let $ \nabla $ be the covariant derivative, then $ \nabla $ carries sections of $B[\lambda]$ to sections of $T^{*}M \otimes B[\lambda]$, so we have
$$\nabla : B[\lambda]{\to} T^{*}M \otimes B[\lambda],$$
where the tensor product $ T^{*}M \otimes B[\lambda]$ can be decomposed under O(n) into irreducible factors. This decomposition defines a natural projection as follows:
$$ B[\lambda]\stackrel{\nabla}{\to}T^{*}M \otimes B[\lambda]\cong\bigoplus_{\beta}B[\beta] \stackrel{Proj}{\to}B[\beta].$$

\textbf{Characteristics of this decomposition}
\begin{itemize}
\item This decomposition $\bigoplus_{\beta}B[\beta]$ is irreducible under the O(n) representation.
\item $\beta $ is the dominant weight.
\item Composition $ G=Proj\circ\nabla  $ is called Generalized gradients or Stein-Weiss operator.
\end{itemize}

\subsection{Young Diagrams}

\begin{figure}[h]
\centerline{\rotatebox{-90}{\epsfxsize =2.5in \epsfysize =3.5in
\epsfbox{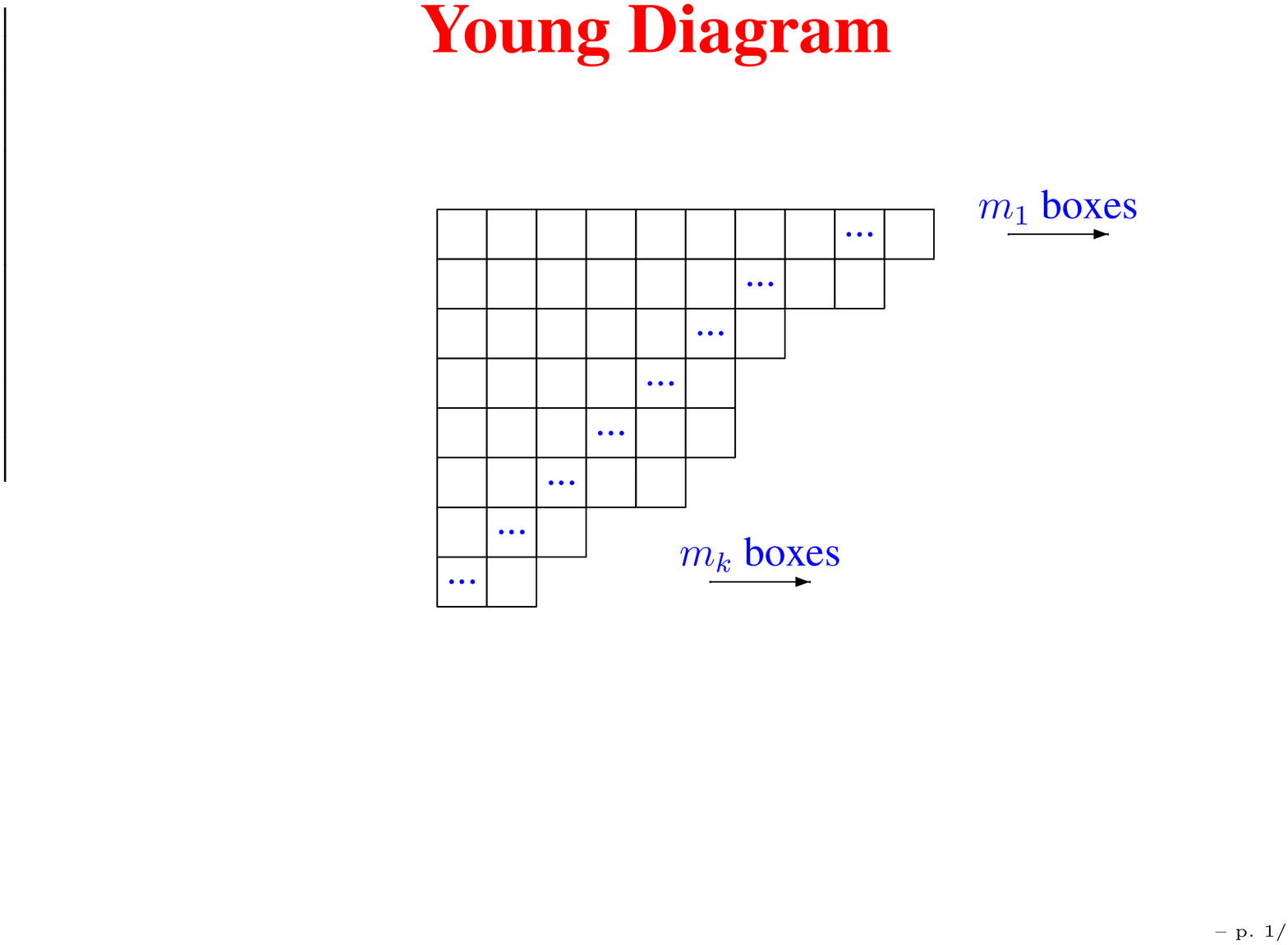}}} \caption{Young Tableaux}
\label{fig1}
\end{figure}

A Young diagram is a representation of partitions of a number n, where partitions are the different ways to express n as sum of natural numbers. It is expressed as a collection of boxes, arranged in left-justified rows, with decreasing number of boxes in each row. For example if $n=m_{1}+m_{2}+...+m_{k},$ is a partition with $ m_{1}\geq m_{2}\geq...\geq m_{k}.$ (See Fig.~\ref{fig1}).

For further lectures in this subject please read \cite{fulton2}. As an application of Young diagrams, in representation theory we can associate a weighted (highest weight) tensor bundle to a Young diagram.

For example, if $ B[2]_{o} $ is a trace-free symmetric  2-tensor bundle, we have the following decomposition:
$$T^{*}M \otimes B[2]_{o} \cong B[3]_{o}\oplus B[2,1]_{o}\oplus B[1].$$
In terms of Young diagrams:
$$T^{*}M \otimes B[2]_{o} \cong \square \square \square_{o}\oplus\frac{\square\square}{\square_{o}} \oplus\square.$$
Notice that $T^{*}M =\bigwedge= B[1]=\square.$

Essentially, the tensor product of $T^{*}M =\bigwedge$ with an irreducible representation B, gives us the direct summands; which are irreducible bundles obtained by adding 1 box and subtracting 1 box. See \cite{murn} for further detail.

\section{Conformal Killing Tensors}
A conformal Killing tensor $V_{a}$ is a vector such that
$$\nabla_{(a}V_{b)}=2Cg_{ab},$$
where C is a constant. The physical understanding is that when the metric is dragged along some congruence of curves, then it remains itself modulo some scale factor "C", which could vary in the manifold. When "C" is zero $V^{a}$ is called a Killing vector, and, clearly the metric is left completely invariant as it is dragged along.
Killing vectors are named for a Norwegian mathematician named W. Killing, who first described these notions in 1892.

A generalization of this concept is defined as follows.

\begin{definition}
A Conformal Killing tensor is a symmetric trace-free tensor field $\sigma_{ab..c}$ with order n, such that
$$[\nabla_{(a}\sigma_{bc...d)}]_{o}=0,$$
this notation means; Trace-Free Part of $\nabla_{(a}\sigma_{bc...d)}=0,$ where the parenthesis means that it is symmetrized in those indices.
\end{definition}
\textbf{Examples}\\
\emph{Case n=1}: $$[\nabla_{(a}\sigma_{b)}]_{o}=0,$$ is
equivalent to
$$\frac{1}{2}(\nabla_{a}\sigma_{b}+\nabla_{b}\sigma_{a})-\frac{1}{n}\nabla^{c}\sigma_{c}g_{ab}=0. $$
\emph{Case n=2}:
$$[\nabla_{(a}\sigma_{bc)}]_{o}=0,$$
is equivalent to
$$\frac{1}{3}(\nabla_{c}\sigma_{ab}+\nabla_{b}\sigma_{ac}+\nabla_{a}\sigma_{bc})-\frac{2}{3(n+2)}(\nabla^{d}\sigma_{dc}g_{ab}
+\nabla^{d}\sigma_{db}g_{ac}+\nabla^{d}\sigma_{da}g_{bc})=0.
$$







\section{The system close and Symmetries}

In this section we explain our method through the proof of 2 theorems. The first theorem concerns conformal Killing tensor of order 1, and the second theorem conformal Killing tensor of order 2. In both cases, we first find a closed system. As an application we find symmetry operators for the Laplace and Yamabe operators in the first case and for the Laplacian in the second case.

We now, give the solution (symmetries $LD$ and $\hat{D}L$) to the equation $LD=\hat{D}L$ from section 1, or equivalently $LD-\hat{D}L=0$ for $L=Y$ and $L=\Delta$.

\begin{theorem}\label{theo1}
Let ($M^{n}$,g) be a Riemannian manifold of dimension n$\geq $3, for any function f in $C^{\infty}(M)$ and $\sigma_{a}$ a conformal Killing tensor, let Y be the Yamabe operator;
$$Yf=\Delta f+\frac{n-2}{4(n-1)}Scf,$$
if
$$
\begin{array}{rl}
(Y )(\sigma^{b}\nabla_{a} + A\nabla_{a}\sigma^{a})f=(\sigma^{b}\nabla_{a} + B\nabla_{a}\sigma^{a})(Y)f,
\end{array}
$$
then
$$A=\frac{n-2}{2n} \qquad and \qquad B =\frac{n+2}{2n}.$$
\end{theorem}

\begin{theorem}\label{theo2}
Let ($M^{n}$,g) be a Riemannian manifold for n$\geq $6, for any function f in $C^{\infty}(M)$ and $\sigma_{ab}$ conformal Killing tensor, then
$$
\begin{array}{rl}
(\Delta)(\sigma^{ab}\nabla_{a}\nabla_{b}+A_{1}\nabla_{a}\sigma^{ab}\nabla_{b}+B_{1}\nabla_{a}\nabla_{b}\sigma^{ab})f=\\\\
(\sigma^{ab}\nabla_{a}\nabla_{b}+A_{2}\nabla_{a}\sigma^{ab}\nabla_{b}+B_{2}\nabla_{a}\nabla_{b}\sigma^{ab})(\Delta)f,
\end{array}
$$
where
$$\begin{array}{rl}
A_{1}=&\frac{n+4}{n+2} \qquad A_{2}=\frac{n}{n+2} \\\\
B_{1}=&\frac{n+4}{4(n+1)}\qquad B_{2}=\frac{n(n-2)}{4(n+1)(n+2)}.\\\\
\end{array}
$$
These values $ A_{1}, A_{2}, B_{1}, B_{2}$ are unique.
\end{theorem}

\begin{pf*}{Proof of theorem ~\ref{theo1}}

A conformal Killing tensor of order 1 satisfies the following overdetermined system:
$$\frac{1}{2}(\nabla_{a}\sigma_{b}+\nabla_{b}\sigma_{a})-\frac{1}{n}\nabla^{c}\sigma_{c}g_{ab}=0,$$
where $ \sigma_{a}$ is a 1-form.
Since $$ \nabla : \Lambda{\to}T^{*}M \otimes \Lambda$$
also by chapter 2 we have,
$$T^{*}M \otimes \Lambda \cong B[2]_{\circ}\oplus \Lambda^{2} \oplus\Lambda^{0}. $$
Then we have three projections, one for each summand. We can repeat this process for each summand as follows:
$$ \nabla : B[2]_{\circ}{\to}T^{*}M \otimes B[2]_{\circ}$$
and
$$T^{*}M \otimes B[2]_{\circ}\cong B[3]_{\circ}\oplus B[2,1]_{\circ} \oplus \Lambda,$$
and we end up with a tree as in fig.~\ref{fig2}.

\begin{figure}[H]
\centerline{\rotatebox{-90}{\epsfxsize =2.5in \epsfysize = 3.5in
\epsfbox{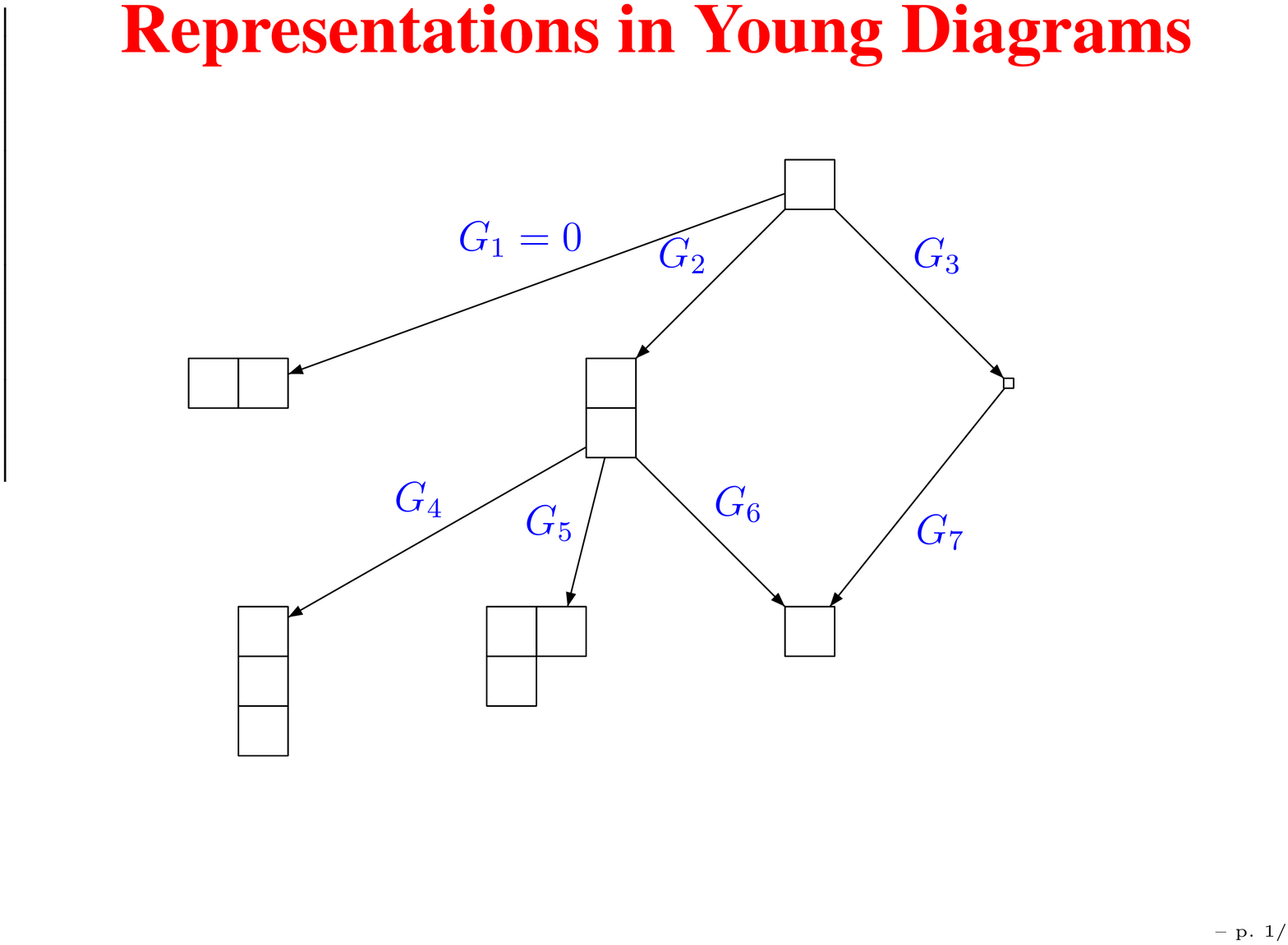}}} \caption{Young Diagram 1}
\label{fig2}
\end{figure}

The gradients $G_{1}$, $G_{2}$, $G_{3}$, are easily computed from the fact that $\nabla = G_{1}+G_{2}+G_{3}$ then
$$
\begin{array}{rl}
(G_{1}\sigma_{d})_{ab}=&\frac{1}{2}(\nabla_{a}\sigma_{b}+\nabla_{b}\sigma_{a})-\frac{1}{n}\nabla^{c} \sigma_{c}g_{ab}\\\\
(G_{2}\sigma_{d})_{ab}=&\frac{1}{2}(\nabla_{a}\sigma_{b}-\nabla_{b}\sigma_{a})=\varphi_{ab}\\\\
(G_{3}\sigma_{d})_{ab}=&\frac{1}{n}\nabla^{c} \sigma_{c} g_{ab}.
\end{array}
$$
Notice that, $(G_{1}\sigma_{d})_{ab}=0$. In general, to find any generalized gradient we use Young symmetrizer, and the projection condition to normalize them. See appendix.

\subsubsection{Strategies Eliminating Higher Derivatives}

We analyze each branch of the tree, and we see that some of them are traces and some are combinations of other branches.\\

\textbf{Composition Type I}: This composition gives us trace terms, so it is zero modulo traces and lower order terms. Then there are no further branches from here.
\begin{figure}[h]
\centerline{\rotatebox{-90}{\epsfxsize =2.5in \epsfysize = 3.5in
\epsfbox{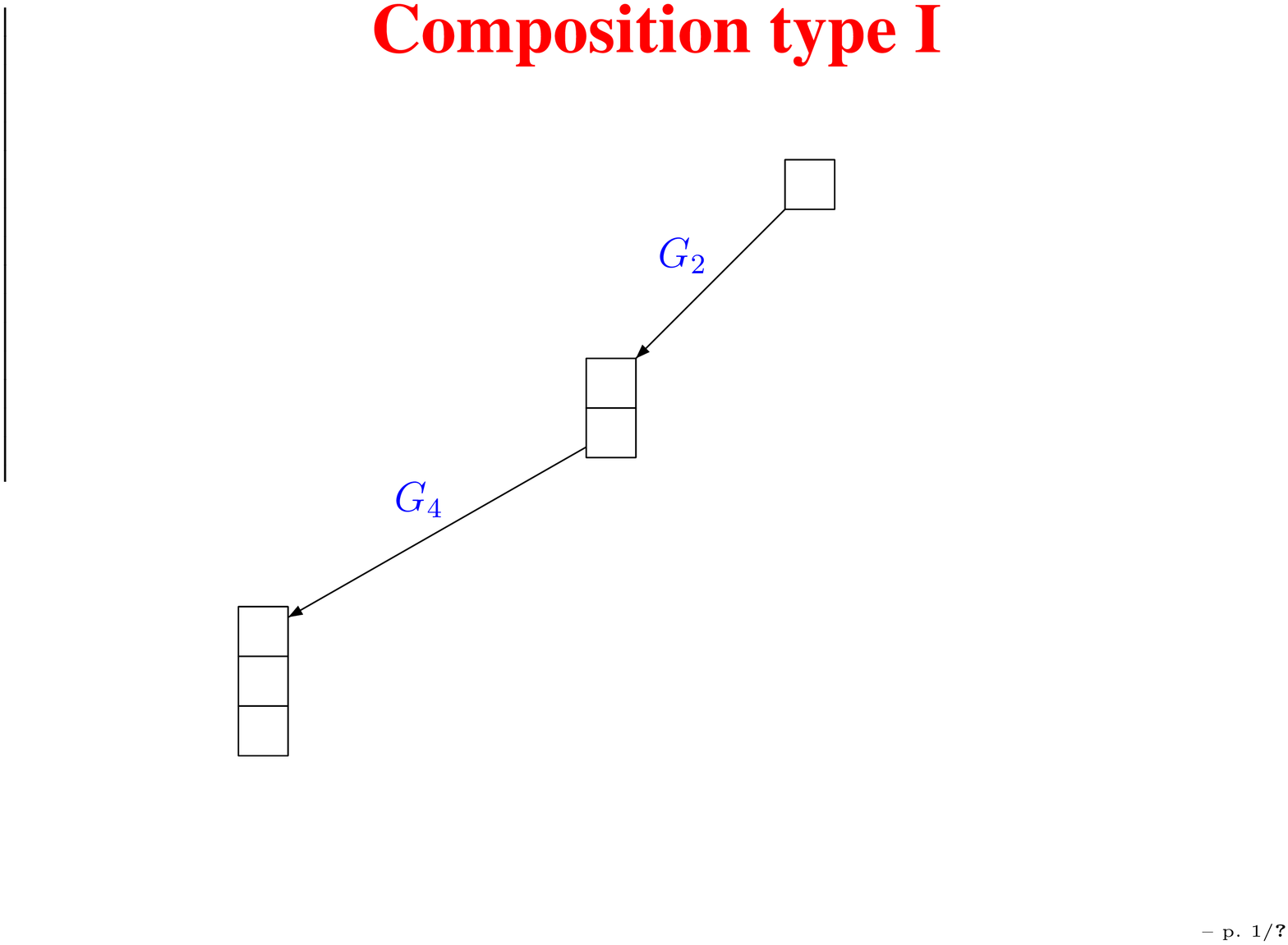}}} \caption{Composition Type I} \label{fig3}
\end{figure}\\

\textbf{Composition Type II}: This allows us identify one branch with another one through composition. In other words one side of the composition can be written in terms of the other side of the composition (see Fig.~\ref{fig4}).
\begin{figure}[h]
\centerline{\rotatebox{-90}{\epsfxsize =2.1in \epsfysize = 3.1in
\epsfbox{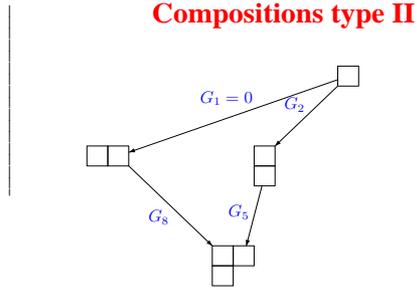}}} \caption{Composition Type II} \label{fig4}
\end{figure}

Iterations of the first and second type compositions in all possible branches determine that all higher derivatives can be written in terms of 4 independent jets, which are $\sigma_{b}$ $\in$ $\Lambda$, $\varphi_{ab}$ $\in$ $\Lambda^{2}$, $\psi$ $\in$ $\Lambda^{2}$ and $\theta_{a}$ $\in$ $\Lambda.$ This is summarized in the next lemma.
\begin{lemma}\label{lem1}
Under hypothesis from theorem ~\ref{theo1} we have the following relations:
$$
\begin{array}{rl}
\nabla^{a}\sigma^{b}=&\frac{1}{2}(\nabla^{a}\sigma^{b}+\nabla^{b}\sigma^{a})-\frac{1}{n}\nabla_{c}\sigma^{c}g^{ab}
+\frac{1}{2}(\nabla^{a}\sigma^{b}-\nabla^{b}\sigma^{a})+\frac{1}{n}\nabla_{c}\sigma^{c}g^{ab}\\\\
=&\frac{1}{2}\varphi^{ab}+\frac{1}{n}\psi g^{ab}\\\\
\nabla^{a}\psi=&\theta^{a}\\\\
\nabla^{a}\varphi^{bc}=&0\\\\
\nabla_{a}\theta^{b}=&-\frac{n}{2(n-1)}Sc^{b}\sigma_{b}-\frac{1}{n-1}\psi Sc
\end{array}
$$
where $\sigma_{b}$ $\in$ $\Lambda$, $\varphi_{ab}$ $\in$  $\Lambda^{2} $, $\psi$ $\in$  $\Lambda^{2} $, and $\theta_{a}$ $\in$ $ \Lambda$ (see Fig. ~\ref{fig5}).
\end{lemma}
\begin{figure}[h]
\centerline{\rotatebox{-90}{\epsfxsize =2.2in \epsfysize = 3.2in
\epsfbox{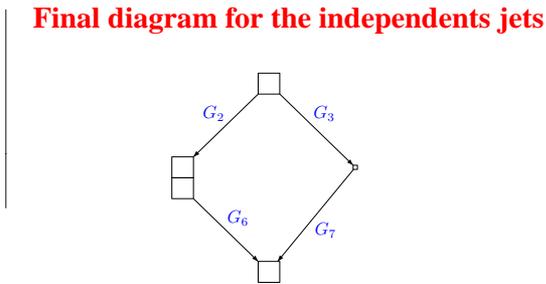}}} \caption{Young Diagram 2}
\label{fig5}
\end{figure}
\subsubsection{Finding symmetries}
To find first order higher symmetry operators of the Yamabe operator, we first set the following equation: $$YDf-\widehat{D}Yf=0,$$ where Y, D and $\widehat{D}$ are:
$$
\begin{array}{rl}
Yf=&\nabla^{a}\nabla_{a}f + \frac{n-2}{4(n-1)}Sc f\\\\
Df=&\sigma^{b}\nabla_{a}f + A_{1}\nabla_{a}\sigma^{a}f\\\\
\widehat{D}f=&\sigma^{b}\nabla_{a}f + A_{2}\nabla_{a}\sigma^{a}f
\end{array}
$$
After expressing all derivatives in term of the 4 independent jets found in lemma ~\ref{lem1} we get:
$$
\begin{array}{rl}
0=YDf-\widehat{D}Yf=
&\frac{2-n+2A_{1}n}{4(n-1)}\sigma_{a}Sc^{a}f+\frac{2-n+2A_{1}n}{4(n-1)}\theta_{a}\nabla_{a}f\\\\
&-\nabla_{a}\nabla_{b}f\varphi_{ab}-\frac{2+A_{1}n-A_{2}n}{n}\psi\nabla^{a}\nabla_{a}f\\\\
&+\frac{2A_{1}+2A_{2}+A_{1}n-A_{2}n}{4(n-1)}\psi Sc f,
\end{array}
$$
so
$$A_{1}=\frac{n-2}{2n} \qquad A_{2}=\frac{n+2}{2n}.$$
Notice that this holds for the curved and non-curved cases. This proves theorem ~\ref{theo1} $ \Box$
\end{pf*}

\begin{pf*}{Proof of theorem ~\ref{theo2}}
Let $ \sigma_{ab}$ be a Killing tensor, then,
$$\frac{1}{3}(\nabla_{c}\sigma_{ab}+\nabla_{b}\sigma_{ac}+\nabla_{a}\sigma_{bc})
-\frac{2}{3(n+2)}(\nabla^{d}\sigma_{dc}g_{ab}+\nabla^{d}\sigma_{db}g_{ac}+\nabla^{d}\sigma_{da}g_{bc})=0
$$
or
\begin{equation}\label{eq1}
[\nabla_{(a}\sigma_{bc)}]_{\circ}=0
\end{equation}
First of all, we want to write all jets of $ \sigma$ in terms of the finitely many jets
\begin{equation}\label{eq2}
\begin{array}{ll}
\qquad \qquad \qquad \qquad \qquad  \sigma_{ab} \in B[2]_{\circ} \\
\qquad \qquad \mu_{cab}\in B[2,1]_{\circ} \qquad \qquad \nu_{a} \in B[1]=\Lambda^{1} \\
\rho_{abcd}\in B[2,2]_{\circ} \qquad \theta_{ab}\in B[2]_{\circ} \qquad \omega_{ab}\in B[1,1] \qquad \varphi \in \bigwedge^{o}\\
\qquad \qquad \beta_{abc} \in B[2,1]_{\circ} \qquad \qquad \lambda_{a} \in B[1] \\
\qquad \qquad \qquad \qquad \qquad  \tau_{ab} \in B[2]_{\circ}, \\
\end{array}
\end{equation}
here $B[2]_{\circ}$ is the bundle of symmetric 2-tensors, $B[2,2]_{\circ}$ is the bundle of algebraic Weyl tensors in symmetrized (as opposed
to antisymmetrized) form. A section is totally trace free and satisfies
$$ \xi_{abcd}= \xi_{(ab)cd}=\xi_{cdab} \qquad \xi_{(abc)d}= 0.$$

The level in the array (~\ref{eq2}) indicates the number of derivatives of $\sigma $ involved. In fact,
\begin{equation}\label{eq3}
\begin{array}{ll}
\qquad \qquad \qquad \qquad \qquad \qquad  \sigma_{ab}\\
\qquad \qquad \mu_{cab}=G_{[2,1][2]}\sigma_{ab} \qquad \qquad \nu_{a}=G_{[1][2]}\sigma_{ab} \\
\rho_{abcd}=G_{[2,2][2,1]}u_{cab}\qquad\theta_{ab}=G_{[2][1]}\nu_{a}\qquad\omega_{ab}=d\nu_{a}\qquad\varphi=\delta\nu_{a} \\
\qquad \qquad \beta_{abc}=G_{[2,1][1,1]}\omega \qquad \qquad \lambda_{a}=d \varphi \\
\qquad \qquad \qquad \qquad \qquad  \tau_{ab} =G_{[2][1]}\lambda.
\end{array}
\end{equation}
Since $\nabla $ carries sections bundles of $B[2]_{\circ}$ to
$$T^{*}M\otimes B[2]_{\circ}\cong B[3]_{\circ}\oplus B[2,1]_{\circ}\oplus \Lambda,$$
and by (~\ref{eq1}), the second level of (~\ref{eq2}) holds the independent first derivatives of $\sigma.$\\
The work begins at the next level. All first derivatives of $\nu$ are recorded in $\theta, \omega, \varphi.$ There are 5 bundles holding the derivatives of $\mu$:
\begin{equation}\label{eq4}
T^{*}M\otimes B[2,1]_{\circ} \cong B[2,2]_{\circ}\oplus B[2,1,1]_{\circ} \oplus B[3,1]_{\circ}\oplus B[2]_{\circ}\oplus B[1,1].
\end{equation}
The part valued in $B[3,1]_{\circ}$ is $G_{[3,1][2,1]}\mu=G_{[3,1][2,1]}\circ G_{[2,1][2]}\sigma.$\\

But module order zero terms, there is one linear relation between the two composition in the diagram

\begin{equation}\label{eq5}
\begin{array}{c}  \\ G_{[3][2]}\\ \
\end{array}
\begin{array}{c} B[3]_{\circ} \\ \uparrow \\ B[2]_{\circ}
\end{array}
\begin{array}{c} \stackrel{G_{[3,1][3]}}{\to} \\ \\ \stackrel{G_{[2,1][2]}}{\to}
\end{array}
\begin{array}{c} B[3,1]_{\circ} \\ \uparrow \\ B[2,1]_{\circ}
\end{array}
\begin{array}{c}  \\ G_{[3,1][2,1]} \\ \
\end{array}
\end{equation}
That is
$$G_{[3,1][2,1]}\circ G_{[2,1][2]}\sigma \sim G_{[3,1][3]}\circ G_{[3][2]}\sigma \sim 0,$$
since $G_{[3][2]}\sigma=[\nabla_{(a}\sigma_{bc)}]_{\circ}=0.$ See Fig.~\ref{fig6}.

\begin{figure}[h]
\centerline{\rotatebox{-90}{\epsfxsize =2.3in \epsfysize = 3.3in
\epsfbox{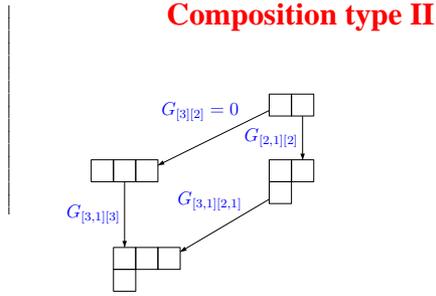}}} \caption{Composition Type II}
\label{fig6}
\end{figure}

The $B[2,2]_{\circ}$ part of $\nabla\mu$ has been given the name $\rho=G_{[2,2][2,1]}\mu$; it is an independent jet.

Next, consider the part of $\nabla\mu$ valued in $B[2,1,1]_{\circ}$;
$$ G_{[2,1,1][2,1]}\mu \sim G_{[2,1,1][2,1]}\circ G_{[2,1][2]}\sigma \sim 0.$$
See Fig.~\ref{fig7}.

\begin{figure}[h!]
\centerline{\rotatebox{-90}{\epsfxsize =2.5in \epsfysize = 3.5in
\epsfbox{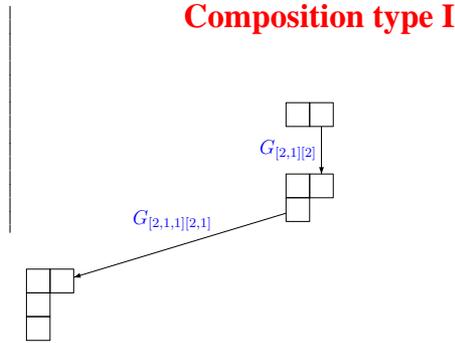}}} \caption{Composition Type I}
\label{fig7}
\end{figure}

The $B[2]_{\circ}$ part of $\nabla\mu$ is an injection into the 4-tensor of the 2-tensor $ (G_{[2,1][2]})^{*}\mu=(G_{[2,1][2]})^{*}\circ G_{[2,1][2]}\sigma$.
By lemma ~\ref{lem2} and $G_{[3][2]}=0$ this is
$$(G_{[2,1][2]})^{*}\circ G_{[2,1][2]}\sigma \sim G_{[2][1]}\circ(G_{[2][1]})^{*}\sigma = G_{[2][1]}\nu=\theta.$$
The $ \Lambda^{2}$ part of $\nabla\mu$ is (carries information equivalent to) another divergent. Let S be
$$ S: \Lambda^{k}\rightarrow B[2,1,...,1];$$
then the $ \Lambda^{2}$ part of $\nabla\mu$ is $ S^{*}\mu=S^{*}G_{[2,1][2]}\sigma.$ Consider the diagram

\begin{equation}\label{eq6}
\begin{array}{c}  \\ G_{[2][1]}^{*}\\ \
\end{array}
\begin{array}{c} \Lambda^{1} \\ \uparrow \\ B[2]_{\circ}
\end{array}
\begin{array}{c} \stackrel{d}{\to} \\ \\
\stackrel{G_{[2,1][2]}}{\to}
\end{array}
\begin{array}{c} \Lambda^{2} \\ \uparrow\\
B[2,1]_{\circ}
\end{array}
\begin{array}{c}  \\ S^{*}\\ \
\end{array}
\end{equation}

The two directions are linearly independent mod terms of order zero, so
$$S^{*}\mu = S^{*}\circ G_{[2,1][2]}\sigma \sim d\circ G_{[2][1]}^{*}\sigma = d\nu = \omega.$$
We move to the next level. $ \nabla\varphi = d\varphi$ is given the name $ \lambda$; this is an independent jet. The three parts of $\nabla\omega$ are $d\omega, \delta \omega $ and $ S\omega$. But $ d\omega = dd\nu=0,$ and $S\omega$ has been given the name $\beta $; this is an independent jet. To handle $\delta \omega,$ we seem to require an essentially third-order argument. Consider third order operator $B[2]_{\circ}\rightarrow \Lambda^{1}.$ Modulo first order operators, there are two linearly independent ones, since there is one map $B[3]_{\circ}\otimes B[2]_{\circ} \rightarrow \Lambda^{1},$ and one map $\Lambda^{1}\otimes B[2]_{\circ} \rightarrow \Lambda^{1}.$ (The symmetric 3-tensor break up into $B[3]_{\circ}$ and $ \Lambda^{1}$ summands.). As a bases of these operators (mod first order) we may take $G_{[2][1]}^{*}\circ G_{[3][2]}^{*}G_{[3][2]}$ and $d\circ G_{[2][1]}^{*}$. Since $G_{[3][2]}$ annihilates $\sigma$,

\begin{equation}\label{eq7}
(\texttt{any} 3^{rd}\texttt{ order op.} B[2]_{\circ} \rightarrow \Lambda^{1})\sigma\ \sim d\delta G_{[3][2]}\circ G_{[2][1]}^{*}\sigma= \lambda.
\end{equation}

In particular,
\begin{equation}\label{eq8}
\delta \omega = \delta d \circ G_{[2][1]}^{*}\sigma \sim \lambda.
\end{equation}

Next, we handle the three parts of $\nabla \theta. $ First,
\begin{equation}\label{eq9}
G_{[3][2]}\theta=G_{[3][2]}\circ G_{[2][1]}\circ(G_{[2][1]})^{*}\sigma \sim b \Delta \sigma \sim \Delta \circ G_{[3][2]} \sigma=0.
\end{equation}
The introduction of the Laplacian $\Delta$ here is via (~\ref{eq1}) together with $G_{[3][2]} \sigma=0$ and lemma ~\ref{lem3} (see Appendix).

The last $ \sim$ in (~\ref{eq9}) is via the "sliding Laplacian principal." The combinatorial formula for the operator $ G_{[p+1][p]}$ involves free indices to the right of the operator, and $\Delta$ has only bound indices.

Next, we handle $G_{[2,1][2]}\theta $ by reversing the vertical arrows in (~\ref{eq6}):

\begin{equation}\label{eq10}
\begin{array}{c}  \\ G_{[2][1]}\\ \
\end{array}
\begin{array}{c} \Lambda^{1} \\ \downarrow\\
B[2]_{\circ}
\end{array}
\begin{array}{c} \stackrel{d}{\to} \\ \\
\stackrel{G_{[2,1][2]}}{\to}
\end{array}
\begin{array}{c} \Lambda^{2} \\ \downarrow\\
B[2,1]_{\circ}
\end{array}
\begin{array}{c}  \\ S\\ \
\end{array}
\end{equation}

This gives
\begin{equation}\label{eq11}
G_{[2,1][2]}\theta = G_{[2,1][2]}\circ G_{[2][1]}\nu \sim S\circ d\nu = S \omega= \beta.
\end{equation}
We handle $b^{*}\theta$ by (~\ref{eq7}).

$\nabla\rho$ has parts in $B[3,2]_{\circ}$, $B[2,2,1]_{\circ}$ and $B[2,1]_{\circ}$. For the $B[3,2]_{\circ}$ part consider the diagram
\begin{equation}\label{eq12}
\begin{array}{c}  \\ G_{[2,2][2,1]}\\ \
\end{array}
\begin{array}{c} B[2,2]_{\circ} \\ \uparrow\\
B[2,1]_{\circ}
\end{array}
\begin{array}{c} \stackrel{G_{[3,2][2,2]}}{\to} \\ \\
\stackrel{G_{[3,1][2,1]}}{\to}
\end{array}
\begin{array}{c} B[3,2]_{\circ} \\ \uparrow\\
B[3,1]_{\circ}
\end{array}
\begin{array}{c}  \\ G_{[3,2][3,1]}\\ \
\end{array}
\end{equation}
This gives
\begin{equation}\label{eq13}
G_{[3,2][2,2]}\rho = G_{[3,2][2,2]}\circ G_{[2,2][2,1]}\mu \sim G_{[3,2][3,1]}\circ G_{[3,1][2,1]}\mu \sim 0,
\end{equation}
since we showed the underbrace assertion above.

For the $B[2,2,1]_{\circ}$, consider the diagram
\begin{equation}\label{eq14}
\begin{array}{c}  \\ G_{[2,2][2,1]}\\ \
\end{array}
\begin{array}{c} B[2,2]_{\circ} \\ \uparrow\\
B[2,1]_{\circ}
\end{array}
\begin{array}{c} \stackrel{G_{[2,2,1][2,2]}}{\to} \\ \\
\stackrel{G_{[2,1,1][2,1]}}{\to}
\end{array}
\begin{array}{c} B[2,2,1]_{\circ} \\ \uparrow\\
B[2,1,1]_{\circ}
\end{array}
\begin{array}{c}  \\ G_{[2,2,1][2,1,1]}\\ \
\end{array}
\end{equation}
This shows that
\begin{equation}\label{eq15}
\begin{array}{rl}
G_{[2,2,1][2,2]}\rho =& G_{[2,2,1][2,2]}\circ G_{[2,2][2,1]}\mu \sim \\
&G_{[2,2,1][2,1,1]}\circ G_{[2,1,1][2,1]}\mu \sim 0,
\end{array}
\end{equation}
where again the underbrace statement is proved above.

For the $B[2,1]_{\circ}$ part of $\nabla \rho,$ namely $G_{[2,2][2,1]}^{*} \rho = G_{[2,2][2,1]}^{*}\circ G_{[2,2][2,1]}\mu,$ we need to get out the 2 bochner formulas on $B[2,1]_{\circ}$, relating the 5 operators $G^{*}G$ (G a natural first order operator targeted at an irreducible bundle) on $B[2,1]_{\circ}$. We have shown in addition that two of these $G^{*}G$, when applied to $\mu$, produce something $\sim 0.$ Thus
\begin{equation}
\begin{array}{rl}\label{eq16}
G_{[2,2][2,1]}^{*} \rho =& G_{[2,2][2,1]}^{*}\circ G_{[2,2][2,1]}\mu \sim \Delta\mu = \Delta \circ G_{[2,1][2]}\sigma\sim \\
&G_{[2,1][2]}\circ \Delta \sigma \sim G_{[2,1][2]}\circ G_{[2][1]}\circ G_{[2][1]}^{*} \sigma= G_{[2,1][2]}\theta \sim \beta
\end{array}
\end{equation}
the last step by (~\ref{eq11}).

We now move to the next level. $\nabla\lambda$ has 3 parts, namely $d\lambda$, $\delta \lambda$ and $G_{[2][1]}\lambda$. But
$$
\begin{array}{rl}
d\lambda=dd\varphi=0\\\\
\delta\lambda \sim \delta\delta\omega
\end{array}
$$
by (~\ref{eq8}). $G_{[2][1]}\lambda$ is given the name $\tau$; it is an independent jet.

Finally on this level, $\nabla\beta$ has 5 parts. For the $B[3,1]_{\circ}$ part, reverse the vertical arrows in (~\ref{eq12})
\begin{equation}\label{eq17}
\begin{array}{c}  \\ G_{[2,2][2,1]}^{*}\\ \
\end{array}
\begin{array}{c} B[2,2]_{\circ} \\ \downarrow\\
B[2,1]_{\circ}
\end{array}
\begin{array}{c} \stackrel{G_{[3,2][2,2]}}{\to} \\ \\
\stackrel{G_{[3,1][2,1]}}{\to}
\end{array}
\begin{array}{c} B[3,2]_{\circ} \\ \downarrow\\
B[3,1]_{\circ}
\end{array}
\begin{array}{c}  \\ G_{[3,2][3,1]}\\ \
\end{array}
\end{equation}
By (3.16),
\begin{equation}\label{eq18}
G_{[3,1][2,1]}\beta \sim G_{[3,1][2,1]}\circ G_{[2,2][2,1]}^{*}\rho \sim G_{[3,2][3,1]}\circ G_{[3,2][2,2]}\rho \sim 0,
\end{equation}
the last step by (~\ref{eq13}).

The $B[2,2]_{\circ}$ part of $\nabla\beta$ is $$G_{[2,2][2,1]}\beta= G_{[2,2][2,1]}\circ S\omega,$$
since there are no invariant maps $Sym^{2}\otimes\bigwedge^{2}\rightarrow B[2,2]_{\circ}.$\\

For the $B[2,1,1]_{\circ}$ part of $\nabla\beta,$ reverse the vertical arrows in (~\ref{eq14}):
\begin{equation}\label{eq19}
\begin{array}{c}  \\ G_{[2,2][2,1]}^{*}\\ \
\end{array}
\begin{array}{c} B[2,2]_{\circ} \\ \downarrow\\
B[2,1]_{\circ}
\end{array}
\begin{array}{c} \stackrel{G_{[2,2,1][2,2]}}{\to} \\ \\
\stackrel{G_{[2,2,1][2,1]}}{\to}
\end{array}
\begin{array}{c} B[2,2,1]_{\circ} \\ \downarrow\\
B[2,1,1]_{\circ}
\end{array}
\begin{array}{c}  \\ G_{[2,2,1][2,1,1]}^{*}\\ \
\end{array}
\end{equation}
This gives $$G_{[2,1,1][2,1]}\beta \sim G_{[2,1,1][2,1]}\circ G_{[2,2][2,1]}^{*}\rho \sim G_{[2,2,1][2,1,1]}^{*}\circ G_{[2,1,1][2,2]}\rho \sim 0,$$
the last step by (~\ref{eq15}).

For the $B[2]_{\circ}$ part of $\nabla\beta,$ we have $$G_{[2,1][2]}^{*}\beta = G_{[2,1][2]}^{*}\circ S\omega.$$
Reversing the horizontal arrows in (~\ref{eq10}), we have
\begin{equation}\label{eq20}
\begin{array}{c}  \\ S\\ \
\end{array}
\begin{array}{c} \Lambda^{2} \\ \downarrow\\
B[2,1]_{\circ}
\end{array}
\begin{array}{c} \stackrel{\delta}{\to} \\ \\
\stackrel{G_{[2,1][2]}^{*}}{\to}
\end{array}
\begin{array}{c} \Lambda^{1} \\ \downarrow\\
B[2]_{\circ}
\end{array}
\begin{array}{c}  \\ G_{[2][1]}\\ \
\end{array}
\end{equation}
Thus
\begin{equation}\label{eq21}
G_{[2,1][2]}^{*}\beta = G_{[2,1][2]}^{*}\circ S \omega \sim G_{[2][1]}\delta\omega \sim G_{[2][1]}\lambda= \tau,
\end{equation}
by (~\ref{eq8}).

The $\bigwedge^{2}$ part of $\nabla\beta$ is
\begin{equation}\label{eq22}
\begin{array}{rl}
S^{*}\beta &\sim S^{*}\circ G_{[2,1][2]}\theta \sim d\circ G_{[2][1]}^{*}\theta \\
&= d\circ G_{[2][1]}^{*}G_{[2][1]}\circ G_{[2][1]}^{*}\sigma \sim d\lambda = dd\varphi=0,
\end{array}
\end{equation}
by (~\ref{eq11}), (~\ref{eq6}) and (~\ref{eq7}).

We may now move to the next level, where we must show that $\nabla\tau \sim 0.$ The three parts of $\nabla\tau$ are $G_{[2][1]}\tau$, $G_{[2,1][2]}\tau$ and $G_{[2][1]}^{*}\tau.$ First,
$$G_{[2][1]}\tau \sim G_{[2][1]}\circ G_{[2,1][2]}\beta$$
by (~\ref{eq21}). Reversing the horizontal arrows in (~\ref{eq5}), we get
\begin{equation}\label{eq23}
\begin{array}{c}  \\ G_{[3,1][2,1]}\\ \
\end{array}
\begin{array}{c} B[3,1]_{\circ} \\ \uparrow\\B[2,1]_{\circ}
\end{array}
\begin{array}{c} \stackrel{G_{[3,1][3]}^{*}}{\to} \\ \\ \stackrel{G_{[2,1][2]}^{*}}{\to}
\end{array}
\begin{array}{c} B[3]_{\circ} \\ \uparrow\\B[2]_{\circ}
\end{array}
\begin{array}{c}  \\ G_{[3][2]}\\ \
\end{array}
\end{equation}
This shows that
$$G_{[3][2]}\tau \sim G_{[3][2]}\circ G_{[2,1][2]}^{*}\beta \sim G_{[3,1][3]}\circ G_{[3,1][2,1]}\beta \sim 0$$
the last step by (~\ref{eq18}).

For the part, by (~\ref{eq21}),$$G_{[2][1]}^{*}\tau \sim G_{[2][1]}^{*}\circ G_{[2,1][2]}^{*}\beta.$$
Reversing the vertical arrows in (~\ref{eq20}), we get
\begin{equation}\label{eq24}
\begin{array}{c}  \\ S^{*}\\ \
\end{array}
\begin{array}{c} \Lambda^{2} \\ \uparrow\\ B[2,1]_{\circ}
\end{array}
\begin{array}{c} \stackrel{\delta}{\to} \\ \\ \stackrel{G_{[2,1][2]}^{*}}{\to}
\end{array}
\begin{array}{c} \Lambda^{1} \\ \uparrow\\ B[2]_{\circ}
\end{array}
\begin{array}{c}  \\ G_{[2][1]}^{*} \\ \
\end{array}
\end{equation}
This shows that
$$G_{[2][1]}^{*}\tau \sim G_{[2][1]}^{*}\circ G_{[2,1][2]}^{*}\beta \sim \delta\circ S^{*}\beta \sim 0$$
the last step by (~\ref{eq22}).

Finally, for the $G_{[2,1][2]}\tau$ part, by (~\ref{eq10}),
$$ G_{[2,1][2]}\tau=G_{[2,1][2]}\circ G_{[2][1]} \lambda \sim S\circ d \lambda = Sdd\varphi=0.$$
So there are not independent jets of $\sigma$ beyond those given in (~\ref{eq2}).

\subsubsection{Finding symmetry}
$$
\begin{array}{rl}
YDf-\widehat{D}Yf=&\frac{2A_{2}-4B_{2}+2A_{2}n-4B_{2}n}{2(n+1)}\lambda^{a}\nabla_{a}f
-A_{2}\nabla_{a}\nabla_{b}f\omega^{ab}\\\\
&+\frac{2A_{2}+B_{1}n-B_{2}n}{n}\varphi \nabla^{a}\nabla_{a}f+\frac{2A_{2}-n+A_{2}n}{n+2}\theta^{ab}\nabla_{a}\nabla_{b}f\\\\
&+\frac{4-2A_{1}+2A_{2}-2A_{1}n+A_{2}n}{n+2}\lambda^{a}\nabla^{b}\nabla_{b}\nabla_{a}f,
\end{array}
$$
where
$$
\begin{array}{rl}
Df=\sigma^{ab}\nabla_{a}\nabla_{b}f + A_{1}\nabla_{a}\sigma^{ab}\nabla_{b}f+B_{1}\nabla_{a}\nabla_{b}\sigma^{ab}f\\\\
Df=\sigma^{ab}\nabla_{a}\nabla_{b}f + A_{2}\nabla_{a}\sigma^{ab}\nabla_{b}f+B_{2}\nabla_{a}\nabla_{b}\sigma^{ab}f.
\end{array}
$$
then
$$\begin{array}{rl}
A_{1}=&\frac{n+4}{n+2} \qquad A_{2}=\frac{n}{n+2}\\\\
B_{1}=&\frac{n+4}{4(n+1)}\qquad B_{2}=\frac{n(n-2)}{4(n+1)(n+2)}.
\end{array}
$$
\end{pf*}

\appendix

\subsection{Generalized gradients (Stein-Weiss operators) from Trace Free Symmetric 2-Tensors}

Since the tensor bundle where the covariant derivative lands can be decomposed into sub-bundles under an irreducible representation of O(n), the projection of the covariant derivative to these sub-bundles are called generalized gradients. Notice that these operators satisfy two conditions:
\begin{equation}
\begin{array}{ll}
1. \nabla = G_{1} + G_{2} +G_{3}\\
2. (G_{1}\sigma)^{2}_{cab} = (G_{1}\sigma)_{cab}.
\end{array}\label{form1}
\end{equation}
The second condition is by linear algebra.

Since $\nabla $ carries sections bundles of $B[2]_{\circ}$ to
$$T^{*}M\otimes B[2]_{\circ}\cong B[3]_{\circ}\oplus B[2,1]_{\circ}\oplus \Lambda,$$
Lets $G_{1}$ , $G_{2}$ and $G_{3}$ be the generalized gradients from $B[2]_{\circ}$, so they are computed as follows:\\\\
\textbf{Projection on B[2,1] Bundle}\\
Non normalized non-Trace term
$$
\begin{array}{rl}
(G_{1}\sigma)_{cab}&=\frac{1}{2}(\nabla_{c}\sigma_{ab}-\nabla_{a}\sigma_{cb})=\frac{1}{2}[\nabla_{c}\sigma_{ab}-
\frac{1}{2}(\nabla_{a}\sigma_{bc}+\nabla_{b}\sigma_{ac})]\\\\
&= \frac{1}{2}(2\nabla_{c}\sigma_{ab}-\nabla_{a}\sigma_{cb}-\nabla_{b}\sigma_{ca})
\end{array}
$$
Normalized non-trace term:
$$(G_{1}\sigma)_{cab}=K(2\nabla_{c}\sigma_{ab}-\nabla_{a}\sigma_{cb}-\nabla_{b}\sigma_{ca})$$
Using the projection condition, we find $K=\frac{1}{3}$, notice that traces
terms has the form: $ \phi^{d}{}_{d*}g_{**}$, for simplicity we will denote $\phi^{d}{}_{da}=\eta_{a}$. So our projection
would be:
$$(G_{1}\sigma)_{cab}=\frac{1}{3}(2\nabla_{c}\sigma_{ab}-\nabla_{a}\sigma_{cb}-\nabla_{b}\sigma_{ca})
+K_{1}(\eta_{a}g_{bc}+\eta_{b}g_{ac})+K_{2}\eta_{c}g_{ab}
$$
Taking traces:
$$
\begin{array}{l}0=\frac{-1}{3}\eta_{c}-\frac{1}{3}\eta_{c}+2K_{1}\eta_{c}+nK_{2}\eta_{c}\qquad( ab-trace)\\\\
0=\frac{2}{3}\eta_{b}-\frac{1}{3}\eta_{b}+K_{1}(\eta_{b}+n\eta_{b})+K_{2}\eta_{b}
\qquad( ac-trace)
\end{array}
$$
System:
$$
\left\{\begin{array}{ll}0=\frac{-2}{3}+2K_{1}+nK_{2},
\\\\0=\frac{1}{3}+K_{1}(1+n)+K_{2}.
\end{array}\right.
$$
Then $$K_{1}=\frac{-1}{3(n-1)},\qquad K_{2}=\frac{2}{3(n-1)},$$
so our projection is:
$$
\begin{array}{rl}
(G_{1}\sigma)_{cab}=&\frac{2}{3}\nabla_{c}\sigma_{ab}-\frac{1}{3}\nabla_{a}\sigma_{cb}-\frac{1}{3}\nabla_{b}\sigma_{ca}
-\frac{1}{3(n-1)}(\eta_{a}g_{bc}+\eta_{b}g_{ac})\\\\
&+\frac{2}{3(n-1)}\eta_{c}g_{ab}.
\end{array}
$$
\textbf{Projection on B[3] Bundle}
$$
\begin{array}{rl}
(G_{3}\sigma)_{cab}=&\frac{1}{3}(\nabla_{c}\sigma_{ab}+\nabla_{b}\sigma_{ac}+\nabla_{a}\sigma_{bc})
+K_{1}\eta_{c}g_{ab}+\\\\
&K_{2}(\eta_{b}g_{ac}+\eta_{a}g_{bc})
\end{array}
$$
Taking traces:
$$
\begin{array}{l}
0=\frac{1}{3}(\eta_{c}+\eta_{c})+nK_{1}\eta_{c}+2K_{2}\eta_{c}\qquad( ab-trace)\\\\
0=\frac{2}{3}\eta_{b}+K_{1}\eta_{b}+K_{2}(n\eta_{b}+\eta_{b})\qquad( ac-trace)\\
\end{array}
$$
System:
$$
\left\{\begin{array}{ll}0=\frac{2}{3}+nK_{1}+2K_{2},
\\\\0=\frac{2}{3}+K_{1}+K_{2}(1+n).
\end{array}\right.
$$
Then $$K_{1}=\frac{-2}{3(n+2)},\qquad K_{2}=\frac{-2}{3(n+2)},$$
so our projection is:
$$
\begin{array}{rl}(G_{3}\sigma)_{cab}=&\frac{1}{3}(\nabla_{c}\sigma_{ab}+\nabla_{b}\sigma_{ac}+\nabla_{a}\sigma_{bc})
+\frac{-2}{3(n+2)}(\eta_{c}g_{ab}+\eta_{b}g_{ac}\\\\
&+\eta_{a}g_{bc}).
\end{array}
$$
\textbf{Projection on B[1] Bundle}
Since adding all the 3 projections we have the identity ($\nabla= G_{1} + G_{2} +G_{3}$), so our last projections is $ G_{2} =\nabla - G_{1} - G_{3}$, then we have
$$
(G_{2}\sigma)_{cab}=\frac{n}{(n+2)(n-1)}(\eta_{b}g_{ac}
+\eta_{a}g_{bc})-\frac{2}{(n+2)(n-1)}\eta_{c}g_{ab}.
$$

\subsection{Bochner Formulas}

\begin{lemma}\label{lem2}
Under conditions of previous section we have the next formula:
$$
\begin{array}{rl}
((2G_{3}^{*}\circ G_{3}-G_{1}^{*} \circ G_{1}-\frac{2n^{2}}{(n+2)(n-1)}G_{[2][1]}\circ G_{[2][1]}^{*})\sigma)_{bc}=\\\\
r^{u}{}_{b}\sigma_{uc} +r^{u}{}_{c}\sigma_{ub}-2R^{u}{}_{c}{}^{a}{}_{b}\sigma_{au}.
\end{array}
$$
\end{lemma}

\begin{proof}
Taking the gradients computed above we have
$$
\begin{array}{rl}
(G_{3}^{*}\circ G_{3}\sigma)_{bc}=&\nabla^{a}(G_{3}\sigma)_{abc}\\\\
&-\frac{1}{3}(\nabla^{a}\nabla_{a}\sigma_{bc}+\nabla^{a}\nabla_{b}\sigma_{ac}+\nabla^{a}\nabla_{c}\sigma_{ab})\\\\
&+\frac{2}{3(n+2)}(\nabla^{a}\eta_{a}g_{bc}+\nabla_{c}\eta_{b}+\nabla_{b}\eta_{c})\\\\
\approx &-\frac{1}{3}\nabla^{a}\nabla_{a}\sigma_{bc}+\frac{2}{3(n+2)}(2\nabla^{a}\eta_{a}g_{bc}-n(\nabla_{c}\eta_{b}
+\nabla_{b}\eta_{c})),
\end{array}
$$
where $\approx$ means module curvature terms $$\frac{1}{3}(-r^{u}{}_{b}\sigma_{uc}-r^{u}{}_{c}\sigma_{ub}+2R^{u}{}_{c}{}^{a}{}_{b}\sigma_{au}).$$
Furthermore,
$$
\begin{array}{rl}
(G_{1}^{*}\circ G_{1}\sigma)_{bc}=&\nabla^{a}(G_{1}\sigma)_{abc}\\\\
&-\frac{2}{3}\nabla^{a}\nabla_{a}\sigma_{bc}+\frac{1}{3}(\nabla^{a}\nabla_{b}\sigma_{ac}
+\nabla^{a}\nabla_{c}\sigma_{ab})\\\\
&-\frac{2}{3(n-1)}\nabla^{a}\eta_{a}g_{bc}+\frac{1}{3(n-1)}(\nabla_{c}\eta_{b}+\nabla_{b}\eta_{c})\\\\
\approx &-\frac{2}{3}\nabla^{a}\nabla_{a}\sigma_{bc}-\frac{1}{3(n-1)}(2\nabla^{a}\eta_{a}g_{bc}-n(\nabla_{c}\eta_{b}
+\nabla_{b}\eta_{c})).
\end{array}
$$
Finally,
$$
\begin{array}{rl}
(G_{2}^{*}\circ G_{2}\sigma)_{bc}=\frac{1}{(n+2)(n-1)}(2\nabla^{a}\eta_{a}g_{bc}-n(\nabla_{c}\eta_{b}+\nabla_{b}\eta_{c})).
\end{array}
$$
Notice that:
$$ \sigma_{ab}\stackrel{G_{[2][1]}^{*}}{\to} -\eta_{a}\stackrel{G_{[2][1]}}{\to} -\frac{1}{2}(\nabla_{a}\eta_{b}+\nabla_{b}\eta_{a}-\frac{2}{n}\nabla^{d}\eta_{d}g_{ab}),$$
so
$$G_{2}^{*}\circ G_{2}=\frac{2n}{(n+2)(n-1)}G_{[2][1]}\circ G_{[2][1]}^{*}.$$
Thus the exact Bochner formula is
$$
\begin{array}{rl}
((2G_{3}^{*}\circ G_{3}-G_{1}^{*} \circ G_{1}-\frac{2n^{2}}{(n+2)(n-1)}G_{[2][1]}\circ G_{[2][1]}^{*})\sigma)_{bc}=\\\\
r^{u}{}_{b}\sigma_{uc} +r^{u}{}_{c}\sigma_{ub}-2R^{u}{}_{c}{}^{a}{}_{b}\sigma_{au}
\end{array}
$$
or
$$
\begin{array}{rl}
((2G_{[3][2]}^{*}\circ G_{[3][2]}-G_{[2,1][2]}^{*} \circ G_{[2,1][2]}-\frac{2n^{2}}{(n+2)(n-1)}G_{[2][1]}\circ G_{[2][1]}^{*})\sigma)_{bc}=\\\\ r^{u}{}_{b}\sigma_{uc} +r^{u}{}_{c}\sigma_{ub}-2R^{u}{}_{c}{}^{a}{}_{b}\sigma_{au}.
\end{array}
$$
\end{proof}
\begin{lemma}\label{lem3}
$$ G_{[3][2]}^{*}\circ G_{[3][2]} +cG_{[2][1]}\circ G_{[2][1]}^{*} + h^{*}h = \nabla^{*}\nabla = \Delta$$ or
$$ G_{3}^{*}\circ G_{3} +cG_{1}\circ G_{1}^{*} + G_{2}^{*}\circ G_{2} = \nabla^{*}\nabla = \Delta.$$
\end{lemma}

\begin{proof}
Immediately from proof of Lemma ~\ref{lem2}.
\end{proof}

\acknowledgements
The second author would like to thank Professor Michael Eastwood and Palle Jorgensen for their valuable conversations.
This work is part of the second author dissertation thesis and was finished after Prof. Branson passed away. The present article is
dedicated to the memory of Thomas Branson.

\end{article}

\begin{thebibliography}{}

\bibitem[\protect\citeauthoryear{Branson, Cap, Eastwood and Gover}{2005}]{b1}
Thomas. Branson, Andreas Cap, Michael Eastwood, and Rod Gover.
\newblock {Prolongations of Geometric Overdetermined Systems}.
\newblock {\em to appear: Int. J Math.}

\bibitem[\protect\citeauthoryear{Branson}{1997}]{b2}
Thomas. Branson.
\newblock {Stein-Weiss Operators and Ellipticity}.
\newblock {\em J. Functional analysis}, 1997.

\bibitem[\protect\citeauthoryear{Branson and Hijazi}{2002}]{bh}
Thomas. Branson and O. Hijazi.
\newblock {Bochner-Weitznbock Formulas Associated with the Rarita-Schwinger Operator}.
\newblock {\em Int. J. of Math.}, Volume 13 N 2 (2002) 137-182.

\bibitem[\protect\citeauthoryear{Cap}{2002}]{cap}
Andreas Cap and A. Rod Gover.
\newblock {Standard Tractors and the Conformal Ambient Metric Construction}.
\newblock {\em preprint}, arXiv:math.DG/0207016 v1.

\bibitem[\protect\citeauthoryear{Carrol}{1997}]{sean}
Sean M. Carrol.
\newblock { Lecture Notes in General Relativity}.
\newblock Institute for Theoretical Physics, University of California, 1997.

\bibitem[\protect\citeauthoryear{Eastwood}{2005}]{east1}
Michael Eastwood.
\newblock {Higher Symmetries of the Laplacian}.
\newblock {\em Annals of Math}, 161 (2005).

\bibitem[\protect\citeauthoryear{Eastwood}{2004}]{east2}
Michael Eastwood.
\newblock {Representations via Overdetermined Systems}.
\newblock {\em Contemporary Mathematics}, (2004).

\bibitem[\protect\citeauthoryear{Fegan}{1976}]{fegan}
H.D. Fegan.
\newblock {Conformally invariant first order differential operators}.
\newblock {\em Quart. J. Math. Oxford}, 27 (1976), 371–378.

\bibitem[\protect\citeauthoryear{Fulton and Harris}{1991}]{fulton}
W. Fulton and J. Harris.
\newblock {Representation Theory, a First Course}.
\newblock {\em Springer}, 1991.

\bibitem[\protect\citeauthoryear{Fulton}{1997}]{fulton2}
W. Fulton.
\newblock {Young Tableaux}.
\newblock {\em London Mathematical Society}, Student Texts 35.

\bibitem[\protect\citeauthoryear{Gover and Silhan}{}]{rod1}
A. Rod Gover and Josef Silhan.
\newblock {The Conformal Killing Equation on Forms-Prolongations and Applications}.

\bibitem[\protect\citeauthoryear{Gover and Silhan}{}]{rod2}
A. Rod Gover and Josef Silhan
\newblock {Invariant Prolongation of the Conformal Killing Equation on Forms}.

\bibitem[\protect\citeauthoryear{Humphreys}{1975}]{hum1}
J.E. Humphreys.
\newblock {Introduction to Lie Algebras and Representation Theory}.
\newblock {\em Grad. Texts Math.}, vol. 21, Springer 1975.

\bibitem[\protect\citeauthoryear{}{1993}]{kolar}
Ivan Kolar, Peter W. Michor, Jan Slovak.
\newblock { Natural Operations in Differential Geometry}.
\newblock {\em Electronic edition; Springer-Verlag}, 1993.

\bibitem[\protect\citeauthoryear{Kress}{1997}]{kress}
J. Kress.
\newblock {Generalised Conformal Killing-Yano tensors:Application to Electrodynamics}
\newblock Ph.D. dissertation thesis, University of Newcastle, Australia, 1997..

\bibitem[\protect\citeauthoryear{Lee}{version 1.52}]{lee}
John M. Lee.
\newblock {Ricci software}.
\newblock http://www.math.washington.edu/$\sim$lee/Ricci/.

\bibitem[\protect\citeauthoryear{Mikhailov}{}]{Mik}
A. Mikhailov.
\newblock {Notes on higher spin symmetries}.
\newblock {\em arXiv: hep-th/0201019}.

\bibitem[\protect\citeauthoryear{Murnaghan}{1963}]{murn}
F. D. Murnaghan.
\newblock {The Theory of Group Representations}.
\newblock New York Dover Publications, 1963.

\bibitem[\protect\citeauthoryear{Orsted}{2000}]{orsted}
Bent Orsted.
\newblock {Generalized Gradients and Poisson Transforms}.
\newblock {\em S´eminaires et Congr`es 4}, SMF 2000.

\bibitem[\protect\citeauthoryear{Penrose}{1984}]{penrose}
R. Penrose, and W Rindler.
\newblock {Spinors and Space-Time}.
\newblock {\em Cambridge University Press}, Volume 1, 1984.

\bibitem[\protect\citeauthoryear{Semmelmann}{2002}]{sem}
Semmelmann U.
\newblock {Conformal Killing Forms on Riemannian Manifolds}.
\newblock {\em preprint}, 2002.

\bibitem[\protect\citeauthoryear{Spencer}{1969}]{spencer}
D.C. Spencer.
\newblock {Overdetermined Systems of Linear Partial Differential Equations}.
\newblock {\em Bull Amer. Math. Soc.}, 75 (1969) 179-239.

\bibitem[\protect\citeauthoryear{Stein and Weiss}{1968}]{sw}
E. Stein and G. Weiss.
\newblock {Generalization of the Cauchy Riemann equations and representation of the rotation group}.
\newblock {\em Amer. J. Math.} 90 (1968) 163-196.

\bibitem[\protect\citeauthoryear{Yamabe}{1960}]{yamabe}
Hidehiko Yamabe.
\newblock {On a Deformation Structures on Compact Manifolds}.
\newblock {\em Osaka Math. J.}, 12 (1960) 21-37.

\bibitem[\protect\citeauthoryear{Yano and Bochner}{1953}]{yano}
K. Yano and S. Bochner.
\newblock {Curvature and Betti Number}.
\newblock {\em Princeton University Press}, 1953.

\end{thebibliography}
\end{document}